\input amstex
\input amsppt.sty
\magnification\magstep1

\def\ni\noindent
\def\sbs{\subset}

\def\as{\operatorname{asdim}}
\def\diam{\operatorname{diam}}

\def\dist{\operatorname{dist}}
\def\Vert{\operatorname{Vert }}
\def\Edge{\operatorname{Edge }}
\def\dim{\operatorname{dim }}

\def\R{\text{\bf R}}

\def\Z{\text{\bf Z}}

\def\G{\Gamma}

\def\sV{\Cal V}
\def\sU{\Cal U}
\def\sW{\Cal W}

\hoffset= 0.0in
\voffset= 0.0in
\hsize=32pc
\vsize=40pc
\baselineskip=24pt
\NoBlackBoxes
\topmatter
\author
 G. Bell and A. Dranishnikov
\endauthor

\title
On asymptotic dimension of groups acting on trees
\endtitle
\abstract
We prove the following.

\proclaim{Theorem} Let $\pi$ be the fundamental group of a finite graph of
groups with finitely generated vertex groups $G_v$  having $\as G_v\le n$ for all
vertices $v$. Then $\as\pi\le n+1$.
\endproclaim

This gives the best possible estimate for the asymptotic dimension of an
HNN extension and the amalgamated product.
\endabstract

\thanks The second author was partially supported by NSF grant DMS-9971709.
\endthanks

\address  University of Louisville, Department of Mathematics, 328 Natural Sciences Building, Louisville, KY 40292, USA
\newline
\newline\indent University of Florida, Department of Mathematics,
P.O.~Box~118105, 358 Little Hall, Gainesville, FL 32611-8105, USA
\endaddress


\email  gcbell01\@gwise.louisville.edu, dranish\@math.ufl.edu
\endemail

\endtopmatter

\document
\head \S0 Introduction \endhead

The notion of the asymptotic dimension $\as X$ of a metric space
$X$ was introduced by Gromov [Gr, pp. 28--30] ($\as_+X$ in his
notation) for studying asymptotic invariants of discrete groups.
He gave three equivalent definitions. The first two definitions
deal with coverings of metric spaces and allow us to use
topological methods.  The third involves mapping into Hilbert
space and allows the use of more analytical techniques. We will
use all three definitions as well as some characterizations of
$\as$ given in [Dr]. Gromov's first definition is as follows: a
metric space $X$ has {\it asymptotic dimension} $\as X\le n$ if
for any (large) number $d$ one can find $n+1$ uniformly bounded
families $\sU^0,\dots,\sU^n$ of $d$-disjoint sets in $X$ such that
the union $\cup_i\sU^i$ is a cover of $X$. Two sets in a metric
space are called $d$-{\it disjoint} if every point in the first is
at least $d$ away from every point in the second. Asymptotic
dimension is a quasi-isometry invariant. Recall that the two
metric spaces $(X,d_X)$ and $(Y,d_Y)$ are said to be {\it
quasi-isometric} if there exist constants $\lambda, \epsilon$ and
$C$ and a map $f:X\to Y$ such that $Y\subset N_C(f(X))$ and
$$\frac1\lambda d_X(z,w)-\epsilon\le d_Y(f(z),f(w))\le \lambda
d_X(z,w)+\epsilon.$$  Here, $N_C(f(X))$ denotes the open
$C$-neighborhood of $f(X).$

A generating set $S$, $S=S^{-1}$, in a group $\Gamma$ defines the
{\it word metric} on $\Gamma$ by the following rule: $d_S(x,y)$ is
the minimal length of a presentation of the element
$x^{-1}y\in\Gamma$ in the alphabet $S$. If a group $\Gamma$ is
finitely generated, then all metrics $d_S$ with finite $S$ are
quasi-isometric. Thus, one can speak about the asymptotic
dimension $\as\Gamma$ of a finitely generated group $\Gamma$.
Gromov proved [Gr, pg. 31] that $\as\Gamma<\infty$ for hyperbolic
groups $\Gamma$.  In [DJ] it was shown that Coxeter groups have
finite asymptotic dimension. In [BD] we proved that the finiteness of
asymptotic dimension is preserved under amalgamated
products and HNN extensions of groups. Results on the finiteness of
asymptotic dimension became quite important because of a theorem of Yu which
states that the Novikov Higher Signature Conjecture holds for
manifolds whose fundamental group has a finite asymptotic
dimension [Yu].

In this paper we prove a finite asymptotic dimensionality theorem
in the case of a general graph of groups.  Moreover, we
significantly improve the existing estimates for the asymptotic
dimension of amalgamated products and HNN extensions. In
particular, we show that  the asymptotic dimension of the
fundamental group of a graph of groups does not exceed the maximum
of the asymptotic dimension of the vertex groups plus one.

We note that asymptotic dimension is a coarse invariant, i.e. it
is an invariant of the {\it coarse category} introduced in [Ro].
We recall that the objects in the coarse category are metric
spaces and morphisms are coarsely proper and coarsely uniform (not
necessarily continuous) maps. A map $f:X\to Y$ between metric
spaces is called {\it coarsely proper} if the preimage
$f^{-1}(B_r(y))$ of every ball in $Y$ is a bounded set in $X$. A
map $f:X\to Y$ is called {\it coarsely uniform} if there is a
function $\rho:\R_+\to\R_+$, tending to infinity, such that
$d_Y(f(x),f(y))\le\rho(d_X(x,y))$ for all $x,y\in X$. Thus,
quasi-isometric metric spaces are coarsely equivalent. We note
that every object in the coarse category is isomorphic to a
discrete metric space, i.e. to a space with a metric taking values
in the nonnegative integers $\Z_+$.

\head \S1 Asymptotic dimension \endhead

We recall that the {\it order (or multiplicity)} of a cover $\sU$
of a space $X$ is the maximum number $n$ such that every $x\in X$
is covered by no more than $n$ elements of $\sU$. If $d>0,$ the
{\it $d$-multiplicity} of a cover $\sU$ of a metric space $X$ is
the maximum number $n$ such that every $d$-ball $B_d(x)$
intersects no more than $n$ elements of $\sU$. Gromov's
second definition of asymptotic dimension follows:

{\it $\as X\le n$ if for any (large) number $d>0$ there is a
uniformly bounded cover $\sV$ of $X$ with $d$-multiplicity $\le
n+1$.}

The equivalence of Gromov's first and second definitions was
proven in [Gr, pp 29--31].

We recall that the Lebesgue number of a cover $\sU$ of a metric
space $X$ is $L(\sU)=\inf\{\max\{d(x,X\setminus U)\mid
U\in\sU\}\mid x\in X\}$. Gromov's second definition can be
reformulated as follows: \proclaim{Assertion 1} $\as X\le n$ if
for any (large) number $d>0$ there is a uniformly bounded cover
$\sU$ of $X$ with multiplicity $\le n+1$ and with Lebesgue number
$L(\sU)> d$.
\endproclaim
\demo{Proof} Let $d$ be given and let $\sV$ be as in Gromov's
second definition. We define $\sU=\{N_d(V)\mid V\in\sV\}$.
Clearly, $L(\sU)\ge d$. Since $B_d(x)$ intersects no more than
$n+1$ elements of $\sV$, no more than $n+1$ elements of $\sU$
contain $x$. Thus, $\sU$ has multiplicity $\le n+1$.

On the other hand, let $d$ be given and let $\sU$ be an open cover
satisfying the conditions of Assertion 1 for $2d$. We define
$\sV=\{U\setminus N_d(X\setminus U)\mid U\in \sU\}$. Since
$L(\sU)\ge 2d$, $\sV$ is still a cover of $X$. The condition
$B_d(x)\cap (U\setminus N_d(X\setminus U))\ne\emptyset$ is
equivalent to $B_d(x) \not\subset N_d(X\setminus U)$. The latter
is equivalent to the condition $x\notin X\setminus U,$ i.e., $x\in
U$. Therefore $\sV$ has multiplicity $\le n+1$. \qed
\enddemo
Let $\sU$ be an open cover of a metric space $X$. The canonical
projection to the nerve $p:X\to Nerve(\sU)$ is defined by the
partition of unity $\{\phi_U:X\to\R\}_{U\in\sU}$, where
$\phi_U(x)=d(x,X\setminus U)/\sum_{V\in\sU} d(x,X\setminus V)$.
The family $\{\phi_U:X\to\R\}_{U\in\sU}$ defines a map $p$ to the
Hilbert space $l_2(\sU)$ with basis indexed by $\sU$. The nerve
$Nerve(\sU)$ of the cover $\sU$ is realized in $l_2(\sU)$ by taking
every vertex $U$ to the corresponding element of the basis.
Clearly, the image of $p$ lies in the nerve.

\proclaim{Proposition 1} For every $k$ and every $\epsilon>0$
there exists a number $\nu=\nu(\epsilon,k)$ such that for every
cover $\sU$ of a metric space $X$ of order $\le k+1$ with Lebesgue
number $L(\sU)>\nu$ the canonical projection to the nerve
$p_{\sU}:X\to Nerve(\sU)$ is $\epsilon$-Lipschitz.
\endproclaim
\demo{Proof}
We take $\nu\ge\frac{(2k+3)^2}{\epsilon}$. Let $x,y\in X$ and $U\in\sU$.
The triangle inequality implies
$$|d(x,X\setminus U)-d(y,X\setminus U)|\le d(x,y).$$
Then
$|\phi_U(x)-\phi_U(y)|\le$
$$\frac{1}{\sum_{V\in\sU}d(x,X\setminus V)}d(x,y)+
d(y,X\setminus U)|\frac{1}{\sum_{V\in\sU}d(x,X\setminus V)}-
\frac{1}{\sum_{V\in\sU}d(y,X\setminus V)}|$$
$$\le\frac{1}{L(\sU)}d(x,y)+
\frac{1}{L(\sU)}(\sum_{V\in\sU}|d(x,X\setminus V)-d(y,X\setminus V)|)\le
\frac{2k+3}{\nu}d(x,y).$$
Then
$\|p(x)-p(y)\|=$
$$=(\sum_{U\in\sU}(\phi_U(x)-\phi_U(y))^2)^{\frac{1}{2}}\le
((2k+2)(\frac{2k+3}{\nu}d(x,y))^2)^{\frac{1}{2}}
\le\frac{(2k+3)^2}{\nu}d(x,y)\le
\epsilon d(x,y).$$\qed
\enddemo

Every countable simplicial complex $K$ can be realized in the
Hilbert space $l_2$ in such a way that all vertices of $K$ are
vectors from some orthonormal basis. As a subspace of the metric
space $l_2,$ $K$ inherits the so-called {\it uniform metric.} All
simplices in this metric are isometric to the standard simplices
$\Delta^n$ (of different dimensions $n$). When the simplicial
complex $K$ is considered with this uniform metric, we denote it
by $K_U.$ Gromov's third definition characterizes asymptotic
dimension by means of Lipschitz maps to uniform polyhedra. We note
that an additional property of the canonical projection is that
the preimage of every (open) simplex $p_{\sU}^{-1}(\sigma)$ lies
in the union of those $U\in\sU$ which correspond to the vertices
of $\sigma$. Thus, if the cover $\sU$ consists of uniformly
bounded sets, then the preimages $p_{\sU}^{-1}(\sigma)$ are also
uniformly bounded. A map $f:X\to K$ to a simplicial complex is
called {\it uniformly cobounded} if the diameters of the preimages
of all simplices $f^{-1}(\sigma)$ are uniformly bounded. Since the
equivalence of Gromov's third definition to other two is not given
in detail in [Gr] we present it here in the form of

\proclaim{Assertion 2} For a metric space $X$ the inequality $\as
X\le n$ is equivalent to the following condition:

For any (small) $\epsilon>0$ there is a uniformly cobounded
$\epsilon$-Lipschitz map $f:X\to K$ to a uniform simplicial
complex of dimension $\le n$.
\endproclaim
\demo{Proof} Let $K$ be an $n$-dimensional simplicial complex. We
consider the cover $\sV_K$ of $K$ by open stars of its vertices
$Ost(v,K)$. Clearly the multiplicity of this cover is $n+1$. If
$K$ is a uniform complex, then there is a uniform lower bound
$c_n$ for the Lebesgue number, i.e., $L(\sV_K)>c_n$ for all $K$
with $\dim K\le n$. For every uniformly cobounded
$\epsilon$-Lipschitz map $f:X\to K$ to a uniform simplicial
complex $K$ of dimension $\le n$ the preimage $f^{-1}(\sV_K)$ will
be a uniformly bounded cover of $X$ with multiplicity $\le n+1$
and with Lebesgue number $>c_n/\epsilon$. By Assertion 1 we have
$\as X\le n$.

On the other hand, let $\as X\le n$. Given $\epsilon>0$ we take
$\nu=\nu(\epsilon, n)$ from Proposition 1. By Gromov's second
definition of $\as$ we can take a uniformly bounded cover $\sU$ of
$X$ of multiplicity $\le n+1$ with Lebesgue number $>\nu$. Then by
Proposition 1 the projection $p_{\sU}:X\to Nerve(\sU)$ is an
$\epsilon$-Lipschitz map to a uniform $n$-dimensional polyhedron.
The preimage $p_{\sU}^{-1}(\sigma)$ of every simplex is contained
in a union of elements of $\sU$ having a common point. Hence
$p_{\sU}$ is uniformly cobounded.\qed
\enddemo

We will use the following theorems proven in
[BD].

 \proclaim{Finite Union Theorem}
Suppose that a metric space is presented as a union $A\cup B$ of subspaces.
Then $\as A\cup B\le\max\{\as A,\as B\}$.
\endproclaim
A family of metric spaces $\{F_{\alpha}\}$ satisfies the
inequality $\as F_{\alpha}\le n$ {\it uniformly} if for arbitrary
large $d>0$ there exist $R$ and $R$-bounded $d$-disjoint families
$\sU^0_{\alpha}\dots\sU^n_{\alpha}$ of subsets of $F_{\alpha}$
such that the union $\cup_i\sU^i_{\alpha}$ is a cover of
$F_{\alpha}$.

\proclaim{Infinite Union Theorem } Assume that
$X=\cup_{\alpha}F_{\alpha}$ and $\as F_{\alpha}\le n$ uniformly.
Suppose that for any $r$ there exists $Y_r\subset X$ with $\as
Y_r\le n$ and such that the family $\{F_{\alpha}\setminus Y_r\}$
is $r$-disjoint. Then $\as X\le n$.
\endproclaim

\head \S2 Simplicial mapping cylinders \endhead

An {\it orientation} on a simplicial complex $K$ is
a linear order on its vertices.
Every subcomplex of an oriented complex has the induced orientation.

On the product of two ordered sets $A\times B$ we define a partial
ordering called the {\it product order} by the rule:
$(a,b)<(a',b')$ if either $a<a'$ and $b\le b'$ or $a\le a'$ and
$b<b'$.

Let $\sigma$ be an oriented $k$-simplex  $x_0,\dots,x_k$ .
The product order on $\{x_0,\dots,x_k\}\times \{0,1\}$ defines a
triangulation on the prism $\sigma\times[0,1]$ by the following rule: every
chain $(x_{i_0},y_{j_0})<\dots< (x_{i_r},y_{j_r})$, $y_j=0,1$,
spans a simplex. We note that the set of vertices for this
triangulation equals $\{x_0,\dots,x_k\}\times\{0,1\}$.

\proclaim{Proposition 2} For every simplicial map $f:X\to Y$
the mapping cylinder $M_{f}$
admits a triangulation with the set of vertices equal to
the disjoint union of vertices of $X$ and $Y$.
\endproclaim
\demo{Proof} We consider an orientation on $X$. For every simplex
$\sigma$ in $X$ we consider the triangulation of the prism
$\sigma\times[0,1]$ defined above. Since the orientations on
simplices $\sigma$  agree with each other, it defines a
triangulation on $X\times[0,1]$. Let $Z=f(X)$. We consider the
image $M$ of the simplicial map $g:X\times[0,1]\to M$ defined by
means of the map $f$ on $X\times\{1\}$. Then $M_{f}$ is a
simplicial complex as it is the union of the simplicial complexes
$M$ and $Y$ along the common subcomplex $Z$. \qed
\enddemo

By the {\it product metric} on the cartesian product $X\times Y$
of two metric spaces $(X,d_X)$ and $(Y,d_Y)$ we mean the metric
$$d((x_1,y_1),(x_2,y_2))=\sqrt{d_X(x_1,x_2)^2+d_Y(y_1,y_2)^2}$$
where $x_1,x_2\in X$ and $y_1,y_2\in Y$.

\proclaim{Proposition 3} For every $n$ there is a constant $c_n$
so that for any simplicial map $g:X\to Y$ of an oriented $n$-dimensional
simplicial complex $X$
the quotient map  $q:X\times[0,1]\to M_g$ of the product
$X\times[0,1]$ to the mapping cylinder
$M_g$ is $c_n$-Lipschitz where $X$ and $M_g$ are given the uniform metrics and
$X\times[0,1]$ has the product metric.
\endproclaim
\demo{Proof}
We note that the quotient map   $q:X\times[0,1]\to M_g$ can be factored
through the uniformization $q':X\times[0,1]\to (X\times[0,1])_U$
where the triangulation on $X\times[0,1]$ is defined by means of
the orientation on $X$ and $(X\times[0,1])_U$ is a uniform complex.
Thus $q=\bar q\circ q'$ where $\bar q:(X\times[0,1])_U\to M_g$ is
a simplicial map between uniform complexes. Since $\bar q$ is
always 1-Lipschitz, it suffices to show that there is a constant $c_n$
so that for any oriented $n$-dimensional complex $X$ the uniformization
$q':X\times[0,1]\to(X\times[0,1])_U$ is $c_n$-Lipschitz.

For each $k>0,$ by a $k$-{\it prism} we mean a subset of $X\times
[0,1]$ of the form $\sigma\times [0,1]$ for some $k-1$-simplex
$\sigma$ of $X.$  First, for each $k$ there exists a constant
$\lambda_k$ so that a uniformization of each $k$-prism is
$\lambda_k$-Lipschitz.

With the $\lambda_k$ as above, take $c_n=\max\{\lambda_k \mid 1\le
k\le 2n+2\}.$ It remains to show that $q'$ is $c_n$-Lipschitz.  To
this end, let $x$ and $y$ in $X\times [0,1]$ be given.  Take a
simplex $\sigma$ in $X$ of minimal dimension with $x\in
\sigma\times [0,1],$ and a simplex $\tau\in X$ of minimal
dimension with $y\in\tau\times [0,1].$

Since the complex $X$ is considered in the uniform metric, its
vertices span a (infinite) simplex $\Delta_0$. The orientation on
$X$ defines an orientation on $\Delta_0$. Therefore the
triangulation on $X\times [0,1]$ is naturally extendable to a
triangulation on $\Delta_0\times [0,1]$ with the same set of
vertices.

Let $\tilde q$ denote the uniformization map of
$\Delta_0^{(2n+1)}\times [0,1]$. Then $\tilde q\mid_{X\times
[0,1]}=q'$.

Let $\rho$ be a simplex in $\Delta_0$ of minimal dimension so that
$\sigma$ and $\tau$ are faces of $\rho.$  Then, $\dim\rho\le
\dim\sigma +\dim\tau +1\le 2n+1.$  Now, $\tilde q$ is
$c_n$-Lipschitz when restricted to the $k$-prism $\rho \times
[0,1]$. Note that $k\le 2n+2.$ Thus, $d(q(x),q(y)) =d(\tilde
q(x),\tilde q(y))\le c_n d(x,y)$. Thus, $q$ is $c_n$-Lipschitz.
\qed
\enddemo

\proclaim{Proposition 4} Let $A\subset W\subset X$ be  subsets in
a geodesic metric space $X$ such that the $r$-neighborhood $N_r(A)$ is
contained in W and let $f:W\to Y$ be a continuous
map to a metric space $Y$. Assume that the restrictions
$f|_{N_r(A)}$ and $f|_{W\setminus N_r(A)}$ are
$\epsilon$-Lipschitz. Then $f$ is $\epsilon$-Lipschitz.
\endproclaim
\demo{Proof}
Let $x,y\in W$  be two points. The only possibility for difficulties with
checking the Lipschitz condition can occur when $x$ and $y$ lie in
different sets $N_r(A)$ and $W\setminus N_r(A).$ Let $x\in N_r(A)$
and $y\in W\setminus N_r(A)$   and let
$\gamma:[0,d(x,y)]\to X$ be a geodesic segment joining them in $X$.
Then there is a $t\in [0,d(x,y)]$ such that
$\gamma(t)\in\partial N_r(A),$  where $\partial N_r(A)=
\{z\in X\mid d(z,A)=r\}.$
Then $d_Y(f(x),f(y))\le d_Y(f(x),f(\gamma(t)))+d_Y(f(\gamma(t)),f(y))\le
\epsilon d_X(x,\gamma(t))+\epsilon d_X(\gamma(t),y)=\epsilon d_X(x,y)$.
\qed
\enddemo
We recall that $L(\sU)$ denotes the Lebesgue number of a cover
$\sU$ of a metric space. If $\sU$ is a cover of a subset $A\subset
X$ of a metric space, we define $L(\sU)=\inf\{\sup\{d(x,X\setminus
U)\mid U\in\sU\}\mid x\in A\}$. Let $b(\sU)$ denote the least
upper bound for the diameters of $U\in\sU$ and let
$\nu(\epsilon,k)$ be the number defined by Proposition 1.

\proclaim{Lemma 1} Let $r>8\epsilon$ and let $\sV$ and $\sU$ be
covers of the $r$-neighborhood $N_r(A)$ of a closed subset $A$ in
a geodesic metric space $X$ by open sets such that both $\sU$ and
$\sV$ have order $\le n+1,$ $Nerve(\sV)$ is orientable,
and $L(\sU)>b(\sV)>L(\sV)\ge
\nu(\epsilon/4c_n,n)$, where $c_n$ is the constant from Proposition 3.
Then there is an $\epsilon$-Lipschitz map
$f:N_r(A)\to M_g$ to the mapping cylinder supplied with the uniform metric
of a simplicial
map $g:Nerve(\sV)\to Nerve(\sU)$ between the nerves  such that
$f$ is uniformly cobounded,
$f|_{\partial N_r(A)}=q(p_{\sV}|_{\partial N_r(A)},0),$ and
$f|_A=p_{\sU}|_A$, where $p_{\sU}:N_r(A)\to Nerve(\sU)$ and
$p_{\sV}:N_r(A)\to Nerve(\sV)$ are the canonical projections to the
nerves.
\endproclaim
\demo{Proof} We recall that the simplicial complexes $Nerve(\sU)$ and
$Nerve(\sV)$ inherit their metrics from the Hilbert spaces $l_2(\sU)$
and $l_2(\sV)$. By Proposition 1 the projections
$p_{\sU}:N_r(A)\to Nerve(\sU)$ and $p_{\sV}:N_r(A)\to Nerve(\sV)$ are
$\epsilon/4c_n$-Lipschitz. Let $t_x=\frac{2d_X(x,A)}{c_nr}$, we define
$$
f(x)=\cases q(p_{\sV}(x),2-t_x) & if\ d(x,A)> c_nr/2\\
t_xgp_{\sV}(x)+(1-t_x)p_{\sU}(x) & \ otherwise.\\
\endcases
$$
In the first case $q^{_-1}f(x)$ lies in $Nerve(\sV)\times[0,1)$.  The map
$q^{-1}f:N_r(A)\to Nerve(\sV)\times[0,1)$ is $\max\{\epsilon\sqrt 2/4c_n,
2\sqrt 2/c_nr\}$-Lipschitz and is therefore $\epsilon/c_n$-Lipschitz.
In the second case the linear combination is taken in the Hilbert
space $l_2(\sU)$. It belongs to $Nerve(\sU)$, since $gp_{\sV}(x)$ and
$p_{\sU}(x)$ lie in the simplex of $Nerve(\sU)$ spanned by all
$U\in\sU$ that contain $x$. Note that the coefficients in the
linear combination are $2/c_nr$-Lipschitz functions. Therefore, they
are $\epsilon/4c_n$-Lipschitz, as $2/r\le\epsilon/4$.

Proposition 3 implies $q$ is $c_n$-Lipschitz.  Hence the
restriction $qq^{-1}f|_{N_{c_nr}(A)\setminus N_{c_nr/2}}$ is
$\epsilon$-Lipschitz. Since the coefficients $t_x$ and $1-t_x$ are
$\epsilon/4c_n$-Lipschitz and $\|gp_{\sV}(x)\|,\|p_{\sU}(x)\|\le
1$, it follows from the Leibnitz rule that the restriction of $f$
to $N_{c_nr/2}(A)$ is also $\epsilon$-Lipschitz. In view of
Proposition 4 we can conclude that $f$ is an $\epsilon$-Lipschitz
map. Note that $f|_{\partial N_r(A)}=p_{\sV}|_{\partial N_r(A)}$
and $f|_A=p_{\sU}|_A$.

Finally, we show $\diam f^{-1}(q(\sigma\times [0,1]))\le 2b(\sU).$
For $x$ and $y$ in $f^{-1}(q(\sigma\times [0,1])),$ there are
three cases to consider, based on $\dist(x,A)$ and $\dist(y,A).$
If both $\dist(x,A)$ and $\dist(y,A)> c_nr/2,$ then
$f(x)=q(p_{\sV}(x),2-t_x)$ and $f(y)=q(p_{\sV}(y),2-t_y),$ so $x$
and $y$ are in $\sigma$ under the projection to the nerve $p_\sV.$
So, there is some set $V\in\sV$ containing $x$ and $y.$  Thus,
$\dist(x,y)\le b(\sV)<b(\sU).$ Next, if $\dist(x,A)<c_nr/2$ and
$\dist(y,A)\ge c_nr/2,$ then $f(x)$ is the same as above, and
$f(y)$ is a linear combination of elements of a simplex in
$Nerve(\sV),$ with coefficients at most 1.  The map $g$ takes an
open set $V\in\sV$ containing $x$ into a set $U\in\sU$ containing
the elements of the linear combination in the expression for
$f(y).$  Thus, $\dist(x,y)\le 2b(\sU).$ Finally, if both elements
lie in the closed $c_nr/2$ neighborhood of $A,$ then both of the
elements lie in a linear combination of elements of a simplex in
$Nerve(\sU),$ so $\dist(x,y)\le 2b(\sU).$ \qed
\enddemo

\head \S3 Groups acting on trees \endhead

A {\it norm} on a group $A$ is a map $\|\ \|:A\to\Z_+$ such that
$\|ab\|\le\|a\|+\|b\|$ , $\|a^{-1}\|=\|a\|$ and $\|x\|=0$ if and
only if $x$ is the unit in $A$. A set of generators $S\subset A$
with $S=S^{-1}$ defines the norm $\|x\|_S$ as the minimal length
of a presentation of $x$ in terms of $S$. A norm on a group
defines
 a left-invariant metric $d$ by $d(x,y)=\|x^{-1}y\|$.
If $G$ is a finitely generated group and $S$ and $S'$ are two finite generating sets,
then the corresponding metrics $d_S$ and $d_{S'}$ define coarsely equivalent
metric spaces $(G,d_S)$ and $(G,d_{S'})$. In particular,
$\as(G,d_S)=\as(G,d_{S'})$, and we can speak about the asymptotic dimension
$\as G$ of a finitely generated group $G$.

 Assume that a group $\Gamma$ acts on a metric space
$X$. For every $R>0$ we define the $R$-{\it stabilizer} $W_R(x_0)$
of the point $x_0\in X$ to be the set of all $g\in\Gamma$ with
$g(x_0)\in B_R(x_0)$. Here $B_R(x)$ denotes the closed ball of
radius $R$ centered at $x$.

\proclaim{Lemma 2}
Assume that a finitely generated group $\Gamma$ acts by isometries on
a tree $X$ with a base point $x_0$.
Suppose that $\as W_R(x_0)\le k$ for all $R$.
Then $\as\Gamma\le k+1$.
\endproclaim

\demo{Proof}   Let $S=S^{-1}$ be a finite generating set for
$\Gamma$. We define a map $\pi:\Gamma\to X$ by the formula
$\pi(g)=g(x_0)$. Using the uniqueness of a geodesic segment
between any two points in $X$ this map can be canonically extended
to a $\Gamma$-equivariant map $\bar\pi:Y\to X$ of the Cayley graph
$Y$ of the group $\Gamma$. Let $\lambda=\max\{d_X(s(x_0),x_0)\mid
s\in S\}$. We show now that $\pi$ is $\lambda$-Lipschitz. Since
the metric $d_S$ on $\Gamma$ is induced from the geodesic metric
on the Cayley graph, it suffices to check that
$d_X(\pi(g),\pi(g'))\le\lambda$ for all $g,g'\in\Gamma$ with
$d_S(g,g')=1$. Without loss of generality we assume that $g'=gs$,
where $s\in S$. Then
$d_X(\pi(g),\pi(g'))=d_X(g(x_0),gs(x_0))=d_X(x_0,s(x_0))\le\lambda$.
Since every edge in the Cayley graph maps via $\bar\pi$ to a geodesic
segment, it follows that $\bar\pi$ is also $\lambda$-Lipschitz.

Observe that $W_R(x_0)=\pi^{-1}(B_R(x_0)),$ $\gamma
B_R(x)=B_R(\gamma(x))$ and $\gamma(\bar\pi^{-1}(B_R(x)))=
\bar\pi^{-1}(B_R(\gamma(x)))$ for all $\gamma\in\Gamma$, $x\in X$
and all $R$. We note that the spaces $W_R(x_0)$ and
$\bar\pi^{-1}(B_R(x_0))$ are in a finite Hausdorff distance from
each other and therefore, coarsely equivalent. We denote the
latter by $\bar W_R(x_0)$.

Let $c_k$ be the Lipschitz constant defined by uniformization of the product
$K^\prime \times [0,1],$ for a $k$-dimensional simplicial complex $K^\prime$
(see Proposition 3).

Given a (small) number $\epsilon >0,$ we will construct an
$\epsilon$-Lipschitz, uniformly cobounded map $\psi: Y\to K$ to a
uniform $k+1$-dimensional simplicial complex.  Let
$\nu=\nu(\epsilon/4c_k,k)$ be from Proposition 1, and take
$r>\max\{\nu, 8/\epsilon\}.$  We need only consider the
$\Gamma$-orbit of $x_0$ in order to get information about
$\Gamma.$ Since $\as X\le 1,$ we have $\as \Gamma x_0\le 1.$ Thus
there is a uniformly bounded cover $\sW$ of $\G x_0$ with
multiplicity $2$ such that the $\lambda r$-enlargement
$\{N_{\lambda r}(W)\}$ is a cover of $\G x_0$ with multiplicity
$2.$ Indeed, using Gromov's first definition of $\as,$ we can find
two families of uniformly bounded, $d$-disjoint sets $\sW_1$ and
$\sW_2$ whose union covers $\G x_0$ with $d\gg \lambda r;$ then,
put $\sW=\sW_1\cup\sW_2.$ For each $W\in\sW,$ let $\bar W$ denote
$\bar\pi^{-1}(W).$ As $\sW$ is uniformly bounded, there is some
$R>0$ so that for each $W\in\sW$ we can find $x_W\in \G x_0$ with
the property that $N_{\lambda r}(W)\subset B_R(x_W).$ As $\bar\pi$
is $\lambda$-Lipschitz, we have $\bar\pi(N_r(\bar W))\subset
N_{\lambda r}(W).$  Thus, $N_r(\bar W)\subset
\bar\pi^{-1}(B_R(x_W))= \bar W_R (x_W).$  Let $\gamma_W\in\Gamma$
have $\pi(\gamma_W)=x_W.$
 We use the condition $\as \bar
W_R(x_0)\le k$ to construct two families of uniformly bounded
covers $\sV$ and $\sU$ of the set $\bar W_R(x_0)$ both with
multiplicity $\le k+1$ such that the Lebesgue number $L(\sV)$ is
greater than $r$ and $L(\sU)>b(\sV)$, where $b(\sV)$ is an upper
bound for diameters of the cover $\sV.$

First, we
construct a family of uniformly cobounded $\epsilon$-Lipschitz maps $\phi_W:
N_r(\bar W)\to K_W$ to uniform $k$-dimensional simplicial complexes as the
canonical projections to the nerve of $\gamma_W\sV$ restricted to
$N_r(\bar W).$

Next, let $W$ and $W^\prime$ be elements of $\sW$ such that $W\cap W^\prime
\neq\emptyset.$  Thus, in the nerve, $Nerve(\sW),$ the sets $W$ and $W^\prime$ define an
edge $e.$  Put $A_e=W\cap W^\prime,$ and let $\bar A_e$ denote $\bar\pi^{-1}(A_e).$
Observe that
$N_r(\bar A_e)\subset N_r(\bar W)$ and $N_r(\bar W^\prime).$  Let $\sV_e$ and
$\sU_e$ denote the restrictions of the covers $\sV$ and $\sU$ to $N_r(\bar A_e).$
We take the triangulation on the mapping cylinder $M_{g_e}$ of a simplicial
map $g_e: Nerve(\sV_e)\to Nerve(\sU_e)$ as in Proposition 2, and
consider $M_{g_e}$ as a uniform simplicial complex.
As $r>8/\epsilon,$ we may apply Lemma 1 to $\bar A_e\subset Y$ and to the
covers $\sV_e$ and $\sU_e$ to obtain an $\epsilon$-Lipschitz map $f_e:
N_r(\bar A_e)\to M_{g_e}.$
We note that the $N_r(\bar A_e)$ are disjoint for distinct edges in the
nerve.  Indeed, $\bar\pi (N_r(\bar A_e))\subset N_{\lambda r}(A_e)\subset
N_{\lambda r}(W)\cap N_{\lambda r}(W^\prime);$ so, if $N_r(\bar A_e)\cap
N_r(A_{e^\prime})\neq\emptyset,$ then the $\lambda r$-enlargement of $\sW$ would
have multiplicity at least 3.  Thus, for each $W\in \sW$ define
$\psi_W:N_r(\bar W)\to K_W\cup_{W\in e} M_{g_e}=L_W$ to the uniform
complex $L_W,$ with mapping cylinders attached as the union of the
map $\phi_W|_{N_r(\bar W)\setminus\cup_{W\in e} N_r(\bar A_e)}$ and
the restrictions $f_e|_{N_r(\bar A_e)\cap N_r(\bar W)},$ for all edges
$e$ in $Nerve(\sW)$ which contain $W$ as a vertex in the nerve $Nerve(\sW).$

We construct $K$ by gluing together the $L_W.$  Clearly, the dimension of $K$
is at most $k+1.$  The maps $\psi_W:Y\to K$ agree on the common parts $\bar A_e$
so they define a map $\psi:Y\to K.$  The map $\psi$ is $\epsilon$-Lipschitz by
Proposition 4, and uniformly cobounded by the Lemma 1.
\qed
\enddemo
In the next section we show (Lemma 3) that if a group $\Gamma$ acts on a tree with
compact quotient and with $\as \Gamma_x\le n$ for all stabilizers, then
$\as W_R(x_0)\le n$ for all $R$. As a corollary of this and Lemma 2 we obtain
the following:
\proclaim{Theorem 1}
Assume that a finitely generated group $\Gamma$ acts on a tree $X$ with a
compact orbit space with the stabilizers $\Gamma_x$ having
$\as \Gamma_x\le n$ for all vertices $x\in X$. Then $\as\Gamma\le n+1$.
\endproclaim

\head \S4 Graph of groups \endhead

We recall some of the basic constructions and known facts about
graphs of groups.  Our development and notation closely follow
that of [S]. Let $Y$ be a non-empty, connected graph.  To each
$P\in \Vert Y,$ associate a group $G_P$ and to each $y\in \Edge
Y,$ a group $G_y=G_{\bar{y}}$ equipped with two injective
homomorphisms, $\phi_y:G_y\rightarrow G_{t(y)}$ and
$\phi_{\bar{y}}: G_{\bar{y}}\rightarrow G_{i(y)}.$

Define the group $F(G,Y)$ to be the group generated by the elements
of the $G_P$ and the elements
$y\in\Edge Y$ subject to the relations:
$$\bar{y}=y^{-1}$$
and
$$y\phi_y(a)y^{-1}=\phi_{\bar{y}}(a),$$
if $y\in \Edge Y$ and $a\in G_y.$ Let $c$ be a path in $Y$
starting at some vertex $P_0.$ Let $y_1,y_2,\ldots ,y_n$ denote
the edges associated to $c,$ where $t(y_i)=P_i.$  Then the length
of $c$ is $n,$ and we write $i(c)=P_0$ and $t(c)=P_n.$  A word of
type $c$ in $F(G,Y)$ is a pair $(c,\mu)$ where $c$ is a path as
above and $\mu$ is a sequence $r_0, r_1,\ldots,r_n,$ where $r_i\in
G_{P_i},$  The associated element of the group $F(G,Y)$ is
$$|c,\mu|=r_0y_1r_1y_2r_2\cdots y_nr_n\in F(G,Y).$$

Serre gives two equivalent definitions of the fundamental group of
the graph of groups $(G,Y).$  For the first description, let $P_0$
denote a fixed vertex and $G_{P_0}$ the associated group. Let
$\pi=\pi_1(G,Y,P_0)$ be the set of elements of $F(G,Y)$ associated
to a path $c$ in Y with $i(c)=t(c)=P_0.$ Obviously, $\pi \subset
F(G,Y)$ is a subgroup.  For the second description, let $T$ be a
maximal subtree of $Y,$ and define $\pi=\pi_1(G,Y,T)$ to be the
quotient of $F(G,Y)$ by the normal subgroup generated by the
elements $t\in\Edge T.$  If $g_y$ denotes the image of $y\in\Edge
Y$ in $\pi,$ then the group $\pi$ is the group generated by the
groups $G_P$ and the elements $g_y$ subject to the relations
$$g_{\bar y}=g_y^{-1},$$
$$g_y\phi_y(a)g_y^{-1}=\phi_{\bar y}(a),$$
and
$$g_t=1$$
where $a\in G_y,$ and $t\in\Edge T.$  So, in particular,
$\phi_t(a)=\phi_{\bar t}(a)$ for all $t\in\Edge T.$  The
equivalence of the descriptions is proven in [S].

\demo{Examples} (1) If $Y$ is the graph with two vertices $P,Q$
and one edge $y,$ then $\pi_1(G,Y,P)=\pi_1(G,Y,Q)=G_P\ast_{G_y}
G_Q,$ the free product of $G_P$ and $G_Q$ amalgamated over $G_y.$

(2) If $Y$ is the graph with one vertex $P$ and one edge $y$ then
$\phi_{\bar y}(G_y)$ is a subgroup of $G_P,$ and $\pi_1(G,Y,P)$ is
precisely the HNN extension of $G_P$ over the subgroup $\phi_{\bar
y}(G_y)$ by means of $\phi_y\phi_{\bar y}^{-1}.$
\enddemo

We now describe the construction of the so-called Bass-Serre tree
$\tilde X$ on which $\pi$ will act by isometries.  Let $T$ be a
maximal tree and let $\pi_P$ denote the canonical image of $G_P$
in $\pi,$ obtained via conjugation by the path $c,$ where $c$ is
the unique path in $T$ from the basepoint $P_0$ to the vertex $P.$
Similarly, let $\pi_y$ denote the image of $\phi_y(G_{t(y)})$ in
$\pi_{t(y)}.$  Then, set

$$\Vert \tilde X=\coprod_{P\in\Vert Y}\pi/\pi_P$$
and
$$\Edge \tilde X=\coprod_{y\in\Edge Y}\pi/\pi_y.$$

For a more explicit description of the edges, observe that the
vertices $x\pi_{i(y)}$ and $xy\pi_{t(y)}$ are connected by an edge
for all $y\in\Edge Y$ and all $x\in\pi.$ Obviously the stabilizer
of the vertices are conjugates of the corresponding vertex groups,
and the stabilizer of the edge connecting $x\pi_{i(y)}$ and
$xy\pi_{t(y)}$ is $xy\pi_yy^{-1}x^{-1},$ a conjugate of the image
of the edge group. This obviously stabilizes the second vertex,
and it stabilizes the first vertex since $y\pi_yy^{-1}= \pi_{\bar
y}\subset \pi_{i(y)}.$ It is known (see [S]) that the action of
left multiplication on $\tilde X$ is isometric.

Now we will assume that the graph $Y$ is finite and that
the groups associated to the edges and
vertices are finitely generated with some fixed set of
generators chosen for each group.  We let
$S$ denote the disjoint union of the generating sets for
the groups, and require that $S=S^{-1}.$  By the norm
$\|x\|$ of an element $x\in G_P$ we mean the minimal
number of generators in the fixed generating set
required to present the element $x.$  We endow
each of the groups $G_P$ with the word metric given by
$\dist (x,y)=\|x^{-1}y\|.$  We extend this metric
to the group $F(G,Y)$ and hence to the subgroup
$\pi_1(G,Y,P_0)$ in the natural way, by adjoining to $S$ the
collection $\{y, y^{-1}\mid y\in \Edge Y\}.$

\proclaim{Proposition 5} Let $Y$ be a non-empty, finite, connected
graph, and $(G,Y)$ the associated graph of finitely generated
groups.  Let $P_0$ be a fixed vertex of $Y,$ then under the action
of $\pi$ on $\tilde X$, the $R$-stabilizer $W_R(1.P_0)$ is
precisely the set of elements of type $c$ in $F(G,Y)$ with
$i(c)=P_0,$ and $l(c)\le R.$
\endproclaim
\demo{Proof}  First, every $x=|c,\mu|\in K$ with $l(c)\le R$ is
 in $W_R(1.P_0)$ since
$\dist_{\tilde X}(xG_{P_0},G_{P_0})$ is the number of edges in the
path $c$ which is at most $R.$  For the reverse inclusion, it
suffices to show that given a reduced word $x=|c,\mu|,$ with
$l(c)=R,$ the length of any coset representative $xG_{P_0}$ is
also $R.$ Consider the word $xg=r_0y_1r_1\cdots y_nr_ng,$ where
$y_i\in \Edge Y,$ $r_i\in G_{P_i}$ and $g\in G_{P_0}.$ The
relations on $F(G,Y)$ allow for only one type of reduction
involving edges to occur.  It occurs in one of the following two
forms:
$$(1)\text{ }y\phi_y(a)=\phi_{\bar y}(a)y,\text{ or}$$
$$(2)\text{ }y\phi_y(a)y^{-1}=\phi_{\bar y}(a),$$
where $a\in G_y$ and $y\in\Edge Y.$ Clearly only a reduction of
type $(2)$ can decrease the (path) length of the word. So, assume
that there are edges $y_{j+1}^{-1}=y_j$ and through a sequence of
type-(1) reductions, $x$ is translated to
$$r_0y_1\cdots y_jr_j\phi_{y_j}(b)y_{j+1}r^\prime_{j+1}\cdots
y_nr_n^\prime.$$ Now, in order for a type-(2) reduction to occur
here, we must have $r_j\phi_{y_j}(b)=\phi_{y_j}(d)$ for some $d\in
G_{y_j}.$  This implies that $r_j=\phi_{y_j}(db^{-1}),$ which
enables a type-(2) reduction of the original word $x.$  As $x$ was
assumed to be reduced, this cannot occur. \qed
\enddemo

\proclaim{Lemma 3} If a group $\Gamma$ acts on a tree with
compact quotient and finitely generated stabilizers satisfying
$\as \Gamma_x\le n$ for all vertices $x$, then $\as W_R(x_0)\le n$ for all $R$.
\endproclaim
\demo{Proof} Consider the subset $K\subset F(G,Y)$ of all words
whose associated path $c$ has $i(c)=P_0.$  This set acts on the
tree by left multiplication, and the fundamental group is a subset
of $K.$  We consider the $R$-stabilizer $W_R(P_0)$ of the fixed
vertex as a subset of $K,$ and show that this set has asymptotic
dimension at most $n$ by induction.  It will follow then that the
stabilizers in the fundamental group will also have asymptotic
dimension at most $n.$

We proceed by induction.  The base case is clear since $W_0$ is
precisely $G_{P_0}$ and by assumption, $\as W_0\le n.$  By
Proposition 5, $W_R$ is precisely the set of all words in $K$ with
associated path $c$ for which $l(c)\le R.$  Let $H_k\subset K$ be
the set of elements whose paths have length precisely $k.$ Then,
in the above notation, $W_R=\cup_{k\le R} H_k,$ so by the Finite
Union Theorem (\S 1), it suffices to prove that $\as H_k\le n$ for
all $k.$  Since $W_0$ and $H_0$ coincide, we proceed to the
induction step on $H_k.$

Observe that $H_k\subset \cup_{y\in\Edge Y} H_{k-1}yG_{t(y)}.$
Since $Y$ was assumed to be a finite graph, this is a finite
union.  So, by the Finite Union Theorem, it suffices to show that
$\as H_k\cap H_{k-1}yG_{t(y)}\le n,$ for some fixed $y\in\Edge Y.$
Next, we let $Y_r=H_{k-1}yN_r(\phi_y(G_y)),$ where
$N_r(\phi_y(G_y))$ is the $r$-neighborhood of $\phi_y(G_y)$ in
$F(G,Y).$  Now, the set $Y_r$ is coarsely equivalent to
$H_{k-1}y\phi_y(G_y),$ and by the relations on $F(G,Y),$ this set
is $H_{k-1}\phi_{\bar y}(G_y)y=H_{k-1}y.$ Finally, $H_{k-1}y$ is
coarsely equivalent to $H_{k-1},$ which, by the inductive
hypothesis, has asymptotic dimension not exceeding $n.$  Thus, we
conclude that $\as Y_r\le n.$ Consider the sets $xyG_{t(y)},$
indexed by those $x\in H_{k-1},$ which do not end with an element
of $\phi_{\bar y}(G_y).$  Notice that this collection does cover
$H_k\cap H_{k-1}yG_{t(y)}$ since we may obtain $x\phi_{\bar
y}(a)yg$ from $xy\phi_y(a)g,$ which is of the required form.  The
map $G_{t(y)}\mapsto xyG_{t(y)}$ is an isometry in the word
metric, so $\as (xyG_{t(y)})\le n$ uniformly.

To apply the Infinite Union Theorem from \S 1,
we need only show that the family $\{H_k\cap H_{k-1}yG_{t(y)}\setminus Y_r\}$
is $r$-disjoint.  To this end,
let $x\neq x^\prime$ and suppose that $z$ and $z^\prime$ are in
$G_{t(y)}\setminus N_r(\phi_y(G_y)).$
Consider $\dist (xyz,x^\prime yz^\prime)=\|z^{-1}y^{-1}x^{-1}x^\prime
yz^\prime\|.$  Write $z=\phi_y(a)s,$
and $z^\prime=\phi_y(a^\prime)s^\prime,$ where $\|s\|$ and $\|s^\prime\|>r.$
Then,
$$\|z^{-1}y^{-1}x^{-1}x^\prime yz^\prime\|$$
$$=\|s^{-1}y^{-1}\phi_{\bar y}(a^{-1})x^{-1}
x^\prime\phi_{\bar y}(a^\prime) ys^\prime\|.$$
Now, a reduction can only occur in the middle, and
if $\phi_{\bar y}(a^{-1})x^{-1}x^\prime\phi_{\bar y}(a^\prime)$
is not in $\phi_{\bar y}(G_y),$ then
$$\|z^{-1}y^{-1}x^{-1}x^\prime yz^\prime\|>\|s\|+\|s^\prime\|>2r.$$
So, suppose that there is a $b\in G_y$ so that $\phi_{\bar y}(b)=
\phi_{\bar y}(a^{-1})x^{-1}
x^\prime\phi_{\bar y}(a^\prime).$  Then, putting $c=ab(a^\prime)^{-1},$
we see that
$x^\prime=x\phi_{\bar y}(c).$  By construction $x^\prime$ cannot end
with a nontrivial
element of $\phi_{\bar y},$ so that $\phi_{\bar y}(c)=e,$
which implies $x^\prime=x,$ a
contradiction.

We conclude that the families are $r$-disjoint, and therefore,
by the Infinite Union Theorem,
$\as H_k<n,$ so that $\as W_R<n$ for every $R.$   \qed
\enddemo

In view of the Bass-Serre structure theorem [S], every group acting
without inversion on a tree with compact quotient is a fundamental
group of a graph of groups, so we may reformulate
Theorem 1 in terms of graphs of groups.

\proclaim{Theorem 1$^\prime$}
Let $(G,Y)$ be a finite graph of groups with finitely generated
vertex groups satisfying the inequality $\as G_v\le n$ for all vertex groups.
Then for the fundamental group we have $\as\pi_1(G,Y,v_0)\le n+1$ for any
vertex $v_0$.
\endproclaim
The estimate of Theorem 1$'$ is exact for HNN extensions. It is
also exact for the general amalgamated product. For example, $\as
A\ast_CB=2$, where $A=B=C=\Z$ and the inclusions $C\to A$ and
$C\to B$ are equal to the multiplication by 2. Indeed, $A\ast_CB$
is isomorphic to the fundamental group of the Klein bottle. Hence
the group $A\ast_CB$ is coarsely isomorphic to the universal cover
$X$ of the Klein bottle which is equal (topologically) to $\R^2$.
One can show that $X$ is coarsely isomorphic to the Euclidean
plane $\R^2$ and then $\as X=\as A\ast_CB=2$. Another possibility
here is to apply Gromov's estimate for the asymptotic dimension of
a uniformly contractible manifold $X$ [Gr, pg. 32]:  for a
uniformly contractible manifold $X$ without boundary, $\as
X=\operatorname{dim} X.$

We believe that the exact formula in the case of the amalgamated
product should be $\as A\ast_CB=\max\{\as A,\as B,\as C+1\}$
where neither of the monomorphisms $C\to A$ and $C\to B$ is an isomorphism.

\enddocument